\magnification=1200
\rightskip=0pt plus1fil
\input amssym.def
\input epsf

\def\sqr#1#2{{\vcenter{\vbox{\hrule height.#2pt
     \hbox{\vrule width.#2pt height#1pt \kern#1pt
     \vrule width.#2pt} \hrule height.#2pt}}}}
\def\square{\ ${\mathchoice\sqr34\sqr34\sqr{2.1}3\sqr{1.5}3}$}

\def \Rthree{{\Bbb R}^3}
\def \BZ {{\Bbb Z}}

\centerline 
{\bf BRUNNIAN BRAIDS AND SOME OF THEIR GENERALIZATIONS}
\bigskip
\centerline {Theodore Stanford \footnote \dag 
{Research supported in part by the Naval Academy Research Council}}

\centerline {Mathematics Department}
\centerline {United States Naval Academy}
\centerline {572C Holloway Road}
\centerline {Annapolis, MD\ \ 21402}
\medskip
\medskip
\centerline {\tt stanford@nadn.navy.mil}
\bigskip

\bigskip
\noindent
{\bf Abstract.}
We use a variation on the commutator collection
process to characterize those pure braids which
become trivial when any one strand is deleted, or,
more generally, those pure braids which become
trivial when all the strands in any one of a list
of sets of strands is deleted.

\bigskip
\noindent
{\bf 0. INTRODUCTION}
\medskip

A Brunnian link is a tame link of $n$ closed
curves in $\Rthree$ such that deleting any one of
the curves results in a trivial link of $n-1$
components.  (See Rolfsen~[11].)  By analogy, we
shall call a braid of $n$ strands a {\it Brunnian
braid} if deleting any one of the strands produces
a trivial braid of $n-1$ strands.  A Brunnian
braid on more than two strands must clearly be a
pure braid, so we confine our attention to the
pure braid group $P_n \subset B_n$.  Though the
idea of ``deleting a strand'' is topological, we
take here a purely algebraic approach.

Brunnian braids were considered by
Levinson~[7,
8]
under the
name ``decomposable braids''.  More generally, 
he considers ``$k$-decomposable braids'', which
become trivial when any $k$ strands are deleted.
In~[7]
he gives a geometric characterization of such braids,
and in~[8]
he gives algebraic characterizations in the cases
$(n,k) = (3,1), (4,1), (4,2)$.  We shall generalize
these results below.

It is not hard to see that the set of $n$-strand
Brunnian braids is a free normal subgroup of $P_n$.
In the Kourovka notebook~[6],
this was called the subgroup of ``smooth'' braids,
and the problem was posed to give a set of free
generators for smooth braids with a given number of
strands. Johnson~[5]
used the Hall commutator
collection process~[4]
to give a set of
generators modulo any group of the lower central
series.  There has appeared in a conference
proceedings~[3]
an abstract of a solution to the
Kourovka problem, after which the problem was
taken from the notebook.  I have been unable to
locate, however, a paper which follows through on
the abstract.  

In this paper, we give sets of generators for the
Brunnian subgroups, the $k$-decomposable
subgroups, and for a somewhat more general class
of subgroups.  Our generating sets are not
minimal, so we do not address the Kourovka
question.  Our method is a finite variation on the
Hall commutator collection process.  A
similar method was used in [10], in order to
show that an element of an arbitrary group $G$ is
``$n$-trivial'' if and only if it is in the $n$th group
of the lower central series of $G$.

We give our characterization of the Brunnian
subgroup of $P_n$ (Corollary~2.3) 
in terms of {\it
monic commutators} (Definition~1.3).  
We set $[x,y] = x^{-1}y^{-1}xy$ for any group elements $x,y$.
Here are some monic commutators:
$[p_{1,3},p_{2,4}] \in P_4$,
$[[p_{1,2},p_{1,3}],p_{1,4}] \in P_4$, and
$[[[p_{1,2},p_{2,3}],p_{3,4}],[p_{4,5},p_{5,6}]]
\in P_6$.  It is not hard to see that these particular
examples are
all Brunnian braids, since deleting any one strand
trivializes at least one entry in the iterated
commutator, which of course trivializes the whole
commutator.  We define (Definition~1.4 and
Proposition~1.6) the {\it support} of a monic
commutator to be the strands whose indices appear
somewhere in the commutator.  Thus the support of
$[p_{1,2},p_{1,4}] \in P_4$ is $\{1,2,4\}$.
Deleting strands $1$,$2$, or $4$ from this
commutator trivializes it, but deleting strand $3$
does nothing.  We characterize the Brunnian
subgroup of $P_n$ as being generated by all monic
commutators whose support is the whole set of
strands $\{1,2,3 \dots n \}$.  As noted below, the
Brunnian subgroup is not finitely-generated (for
$n>2$), so its list of generators necessarily
includes commutators of with an arbitrary number
of brackets.

More generally, consider any finite collection of subsets
$S_1, S_2, \dots S_m \subset \{1,2, \dots n\}$.  Our main
theorem states that the subgroup of braids which become
trivial when the strands in any one $S_i$ are cut is generated
by the set of monic commutators whose support intersects
each $S_i$ nontrivially.  When $S_i = \{i\}$ we get the
characterization of Brunnian braids, and when the $S_i$
consist of all subsets with $k$ elements, we get a
characterization of $k$-decomposable braids.

\medskip
\noindent
{\bf Appreciation.} \ I would like to than Joan Birman, Tony
Gaglione, and Mark Meyerson for helpful 
remarks and conversations.

\bigskip
\noindent
{\bf 1. BASIC IDEAS}
\medskip

A standard reference on braids is 
Birman~[2].
Recall that the pure braid group $P_n$ is
generated by $p_{a,b}$ for $1 \le a<b \le n$,
where $p_{a,b}$ is the braid which links strand
$a$ and strand $b$ in front of the other strands.
Artin's~[1]
semidirect product
decomposition may be used to find a finite
presentation for $P_n$. Here is one version:

\medskip
\item {A.}
$p_{a,b} p_{a,c} p_{b,c} = p_{a,c} p_{b,c} p_{a,b} 
= p_{b,c} p_{a,b} p_{a,c}$\ \ \ 
for all $1 \le a < b < c \le n$
\smallskip
\item {B.}
$p_{a,b} p_{c,d} = p_{c,d} p_{a,b}$\ \ 
and\ \ $p_{a,d} p_{b,c} = p_{b,c} p_{a,d}$\ \ \ 
for all $1 \le a < b < c < d \le n$
\smallskip
\item {C.}
$p_{a,c} p_{b,c}^{-1} p_{b,d} p_{b,c} = 
p_{b,c}^{-1} p_{b,d} p_{b,c} p_{a,c}$
for all \ \ \ $1 \le a < b < c < d \le n$
\medskip

Each of the above relations corresponds to a 
``geometrically obvious'' commutation relation.
For the (B) relations this is immediate, and examples for
the (A) and (C) relations are shown in the figure below.

\vfil
\eject

\centerline {\hbox {
\epsfysize = .7 truein
\epsffile {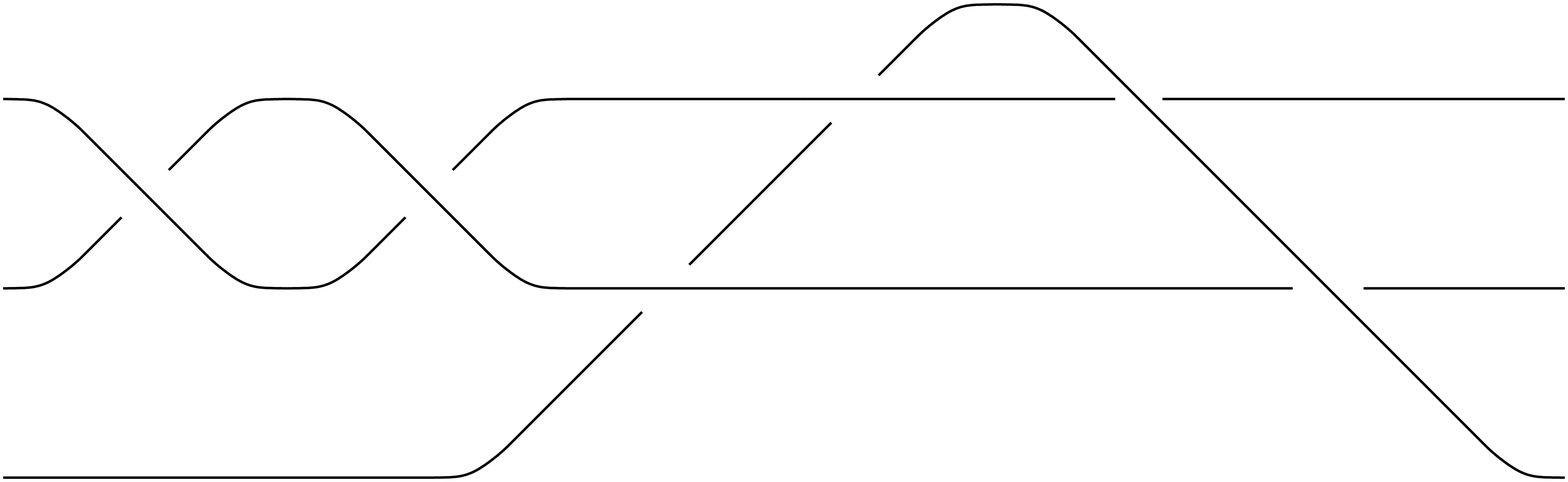}
\hskip .5 truein
\epsfysize = .7 truein
\epsffile {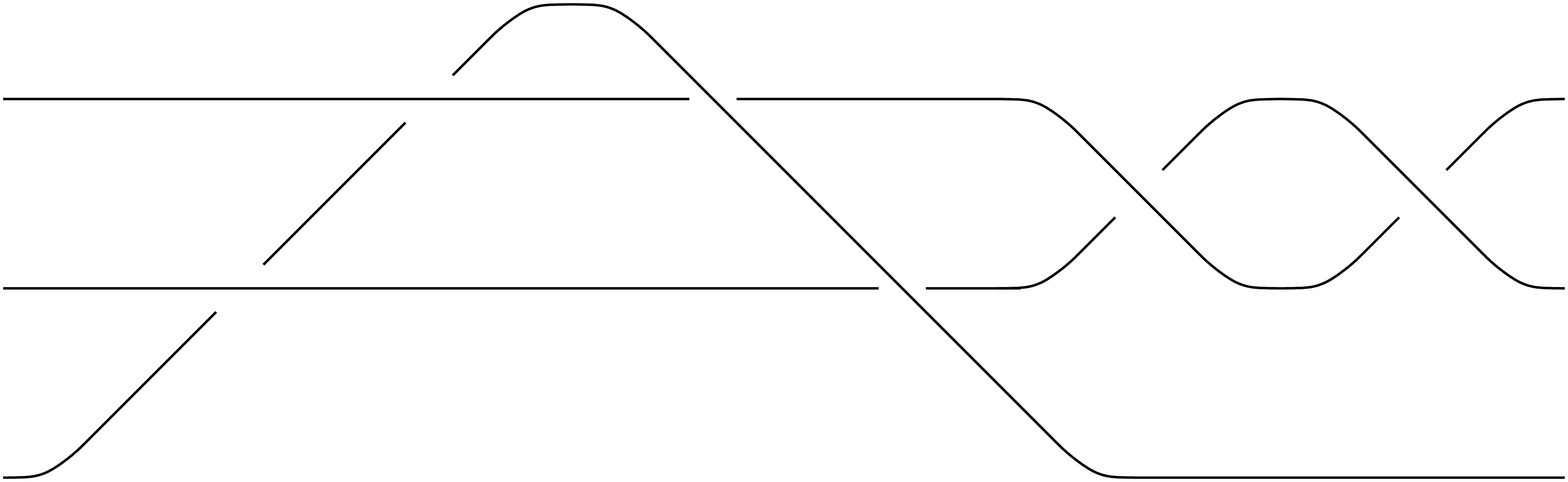}}}
$$p_{1,2} (p_{2,3} p_{1,3}) = (p_{2,3} p_{1,3}) p_{1,2}$$

\medskip
\centerline {\hbox {
\epsfysize = 1 truein
\epsffile {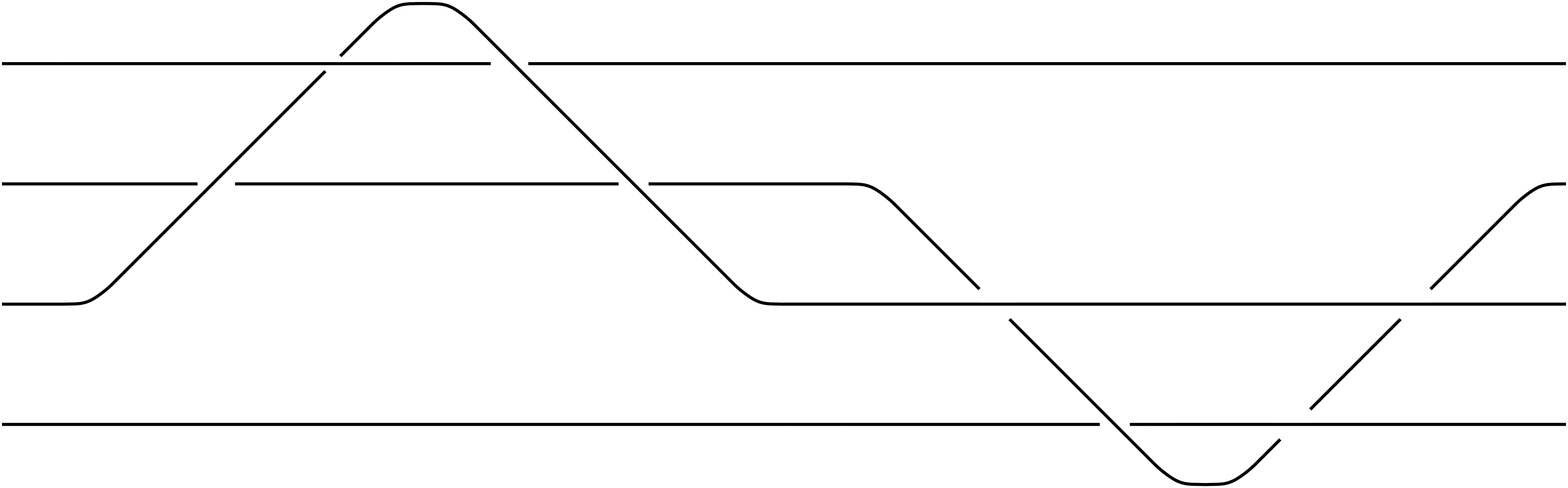}
\hskip .5 truein
\epsfysize = 1 truein
\epsffile {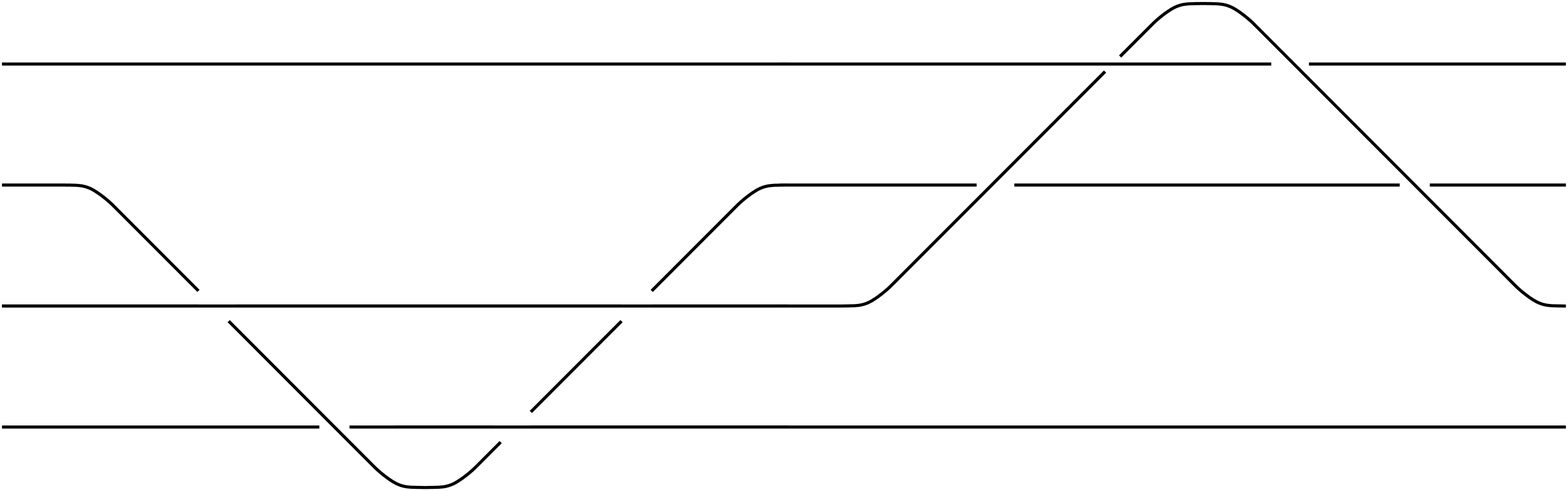}}}
$$p_{1,3} (p_{2,3}^{-1} p_{2,4} p_{2,3}) 
= (p_{2,3}^{-1} p_{2,4} p_{2,3}) p_{1,3}$$
\bigskip

Fix the positive integer $n$, and let 
$N = \{1,2, \dots n\}$. Suppose $S \in n$.  Denote by
$\overline S$ the complement of $S$ in $N$. Let
$P_S$ be subgroup of $P_N = P_n$ generated by $p_{a,b}$
with $a \in S$ and $b \in S$, and let $Q_S$ be the
subgroup generated by $p_{a,b}$ such that either
$a \in S$ or $b \in S$.  The standard semidirect
product decomposition states that there is a
retraction homomorphism from $P_N$ to $P_N$
whose image is $P_{N-1}$ and 
whose kernel is $Q_{\{n\}}$.  This homomorphism is
accomplished geometrically 
by cutting or trivializing
the $n$th string of a pure
braid.  If instead we cut all the strings in
$S \subset N$, then we still get a retraction map and a
semidirect product decomposition.

\bigskip
\noindent
{\bf Proposition 1.1.} \ 
{\it For each $S \subset N$,
there exists a retraction homomorphism $\phi_S: P_n \to P_n$
whose image is $P_{\overline S}$ and whose kernel is $Q_S$.
Moreover, $\phi_{S_1} \circ \phi_{S_2} = \phi_{S_1 \cup S_2}$
for any $S_1, S_2 \subset N$.}

\medskip
\noindent
{\it Proof:} \ 
Define $\phi_S (p_{a,b}) = p_{a,b}$ if $a \notin S$ and 
$b \notin S$,
and $\phi_S (p_{a,b}) = 1$ otherwise.   Then check using
the presentation above that this defines a homomorphism with
the required properties.
\square
\medskip

Now, given a set of subsets $S_1, S_2, S_3, \dots
S_m \subset n$, we are interested
in characterizing $\cap_{i=1}^m Q_{S_i}$.  First
note it follows from Proposition~1.1 that
$P_{S_1} \cap P_{S_2} = P_{S_1 \cap S_2}$.  Inductively, we have

\medskip
\noindent
{\bf Proposition 1.2.}\ \ 
{\it Let $S_1, S_2, \dots
S_m \subset N$.  Then
$\cap_{i=1}^m P_{S_i} = P_{\cap_{i=1}^m S_i}$.}
\medskip

By way of contrast, the subgroup $Q = \cap_{i=1}^m
Q_{S_i}$ is not in general finitely-generated.
For example, consider $Q_{\{2\}}$ and $Q_{\{3\}}$
in $P_3$.  Both of these subgroups are free on two
generators, so $Q_{\{2\}} \cap Q_{\{3\}}$ is also
free.  The retraction $\phi_{\{3\}}: P_3 \to P_3$
restricts to a retraction $Q_{\{2\}}
\to Q_{\{2\}}$, whose kernel is $Q_{\{2\}} \cap Q_{\{3\}}$
and whose image is the infinite cyclic group
$\langle p_{1,2} \rangle$.  
Thus $Q_{\{2\}} \cap Q_{\{3\}}$ has
infinite index in $Q_{\{2\}}$, and therefore
$Q_{\{2\}} \cap Q_{\{3\}}$ is not finitely-generated.

We shall show that $\cap_{i=1}^m Q_{S_i}$ is
generated by a subset of the monic commutators of
$P_n$.

\medskip
\noindent
{\bf Definition 1.3.}\ \ 
A {\it monic commutator} is an element of
$P_n$ defined recursively as follows
\smallskip
\item {A.}
$p_{a,b}$ and $p^{-1}_{a,b}$ are 
monic commutators for all $1 \le a< b \le n$.
\smallskip
\item {B.}
If $x$ and $y$ are monic commutators, 
and $[x,y] \ne 1$, then $[x,y]$ is a monic commutator.

\medskip
\noindent
{\bf Definition 1.4.} \ 
If $x \in P_n$, then
the {\it support} $\sigma (x)$ is 
the intersection of all $S \subset n$
such that $x \in P_S$.
\medskip

We have $\sigma(p_{a,b}^{\pm 1}) = \{a,b\}$.
Two other things are also immediate.  First,
$\sigma ([x,y]) \subset \sigma (x) \cup \sigma (y)$.
Second, if $\sigma (x) \cap S = \emptyset$, then
$\phi_S (x) = x$.  For monic commutators, 
we can say more:

\medskip
\noindent
{\bf Proposition 1.5.}
If $x$ is a monic commutator in $P_n$, then $\phi_S (x) = 1$
if and only $\sigma (x) \cap  S \ne \emptyset$.
\medskip

\noindent
Proof:  
If $\sigma (x) \cap  S = \emptyset$, then 
$\phi_S (x) = x \ne 1$ by definition.
For the converse, it suffices by Proposition~1.1
to show that for a monic commutator $x$, 
$i \in \sigma (x)$ implies that $\phi_{\{i\}} (x) = 1$.
This is certainly true when $x = p_{a,b}^{\pm 1}$.  Suppose that
it is true for two monic commutators $x$ and $y$.
If $i \in \sigma ([x,y])$,
then $i \in \sigma (x)$ or $i \in \sigma (y)$.
Then $\phi_{\{i\}} (x) = 1$ or $\phi_{\{i\}} (y) = 1$,
and in either case $\phi_{\{i\}} ([x,y]) = 1$.
\square

\medskip
\noindent
{\bf Proposition 1.6.} \ 
{\it If $x$, $y$, and $[x,y]$ are monic commutators, then
$\sigma ([x,y]) = \sigma (x) \cup \sigma (y)$.}

\medskip
\noindent
Proof:
If $i \notin \sigma ([x,y])$ then $\phi_{\{i\}} ([x,y]) \ne 1$, 
and therefore $i \notin \sigma (x) \cup \sigma (y)$.
\square
\medskip

\bigskip
\noindent
{\bf 2.  THE THEOREM}

\medskip
\noindent 
{\bf Theorem 2.1.}\ \ 
{\it For $1 \le i \le m$, let $S_i \subset N$.
Then $\cap_{i=1}^m Q_{S_i}$ 
is generated by the set of monic commutators $x$
such that $\sigma (x) \cap S_i \ne \emptyset$
for all $1 \le i \le m$.}
\medskip

We obtain as corollaries characterizations of
$k$-decomposable braids and of Brunnian braids,
generalizing results of 
Levinson~[8].

\medskip
\noindent
{\bf Corollary 2.2.}\ \ 
{\it The normal subgroup of
$k$-decomposable $n$-strand braids is generated in $P_n$ by
all monic commutators whose support has cardinality at least
$n-k+1$.}

\medskip
\noindent
{\bf Corollary 2.3.}\ \ 
{\it The normal subgroup of all Brunnian braids in 
$P_n$ is generated
by all monic commutators whose support is
$N = \{1,2, \dots n\}$.}
\medskip

\medskip
\noindent
{\bf Corollary 2.4.}\ \ 
{\it The subgroup of $(n-2)$-decomposable braids is
the commutator subgroup of $P_n$.}
\medskip

Corollary~2.4 
may also be proved directly by observing that if
$|S| = n-2$ then the image of $\phi_S$ is a two-strand pure
braid group, isomorphic to $\BZ$.  More specifically,
if  $S = N - \{i,j\}$, then $\phi_S (x)$ can be taken
to be the integer which measures the
linking number of strands $i$ and $j$ in the braid $x$.
Then $x$ is in the commutator subgroup of $P_n$ if and 
only if all these linking numbers vanish.

\bigskip
\noindent
{\it Proof of Theorem~2.1:}\ \ 
If $x$ is a monic commutator
which intersects each $S_i$ nontrivially, then $\phi_{S_i}
(x) = 1$ for all $i$, and therefore $x \in \cap_{i=1}^m
Q_{S_i}$.  We need to show that any element in 
$\cap_{i=1}^m Q_{S_i}$
can be written as a product of such commutators.
First we will describe a type of commutator collection
process, where any element of $P_n$ may be written as a
product of monic commutators such that the monic commutators
with common support are grouped together.  When we apply
this to an element of $\cap_{i=1}^m Q_{S_i}$, we will find
that all the monic commutators whose support misses one of
the $S_i$ will drop out.

Let $q \in P_n$ be given as a word in the $p_{i,j}^{\pm 1}$.
Fix a total ordering $T_1, T_2, \dots T_{2^n}$ of the
subsets of $N$, arbitrary except that we require that if
$T_i \subset T_j$ then $i < j$.  For notational convenience
we shall write $T_i < T_j$ if $i<j$.  We claim that there
exist elements $q_1, q_2, \dots q_{2^k} \in P_n$ such that
$q = q_1 q_2 \dots q_{2^k}$, and such that each $q_i$ a
product of monic commutators of support $T_i$,
where an empty product is taken to be $1 \in P_n$.
Suppose inductively that we have written $q = q_1 q_2 \dots
q_r s$, where $q_i$ is a product of monic commutators of
support $T_i$ for $1 \le i \le r$, and $s$ is a product of
monic commutators, each of support $\ge T_r$.  We want to
find all the monic commutators in $s$ of support $T_r$, and
move them back into $q_r$.  We may do this by inserting
commutators, each of which has support $>T_r$.  More
precisely, let $s = x_1 x_2 \dots x_t y z$, where each $x_i$
is a monic commutator of support $>T_r$, $y$ is monic
commutator of of support $T_r$, and $z$ is a product of
monic commutators of support $\ge T_r$.  We may then write
$s = y x_1 [y,x_1] x_2 [y,x_2] x_3 \dots x_t [y,x_t]$, and
then move $y$ into $q_r$.  Even though this increases the
number of monic commutators in $s$, it decreases by one the
number of them that have support $T_r$.  Thus it is possible
to continue until all the monic commutators with support
$T_r$ are contained in $q_r$, and inductively we can
continue until $q = q_1 q_2 \dots q_{2^k}$.

Now let $q \in \cap_{i=1}^m Q_{S_i}$ and $q = q_1 q_2 \dots
q_{2^k}$ as above.  Fix an arbitrary
$i$ between $1$ and $n$.  We need
to show that if $T_j \cap S_i = \emptyset$, then $q_j = 1$.
This is done by induction on $|T_j|$.  Suppose that 
$q_{j^\prime} = 1$ for all $T_{j^\prime} \subset T_j$.
By Proposition~1.1, 
$\phi_{\overline T_j}$ factors through $\phi_{S_i}$,
so $$1 = \phi_{\overline T_j} (q)
= \phi_{\overline T_j} (q_1)
\phi_{\overline T_j} (q_2)
\dots
\phi_{\overline T_j} (q_{2^n})$$
and the only possible nontrivial element of this product is
$\phi_{\overline T_j}(q_j)$, which must then be trivial as well.
\square

\vfil
\eject

\noindent
{\bf References.}

\smallskip
\item {[1]}
E. Artin.
{\it Theory of braids}.
Annals of Mathematics (2) 48 (1947), 101--126.

\smallskip
\item {[2]}
J.S. Birman.
``Braids, Links and Mapping Class Groups.''
Annals of Mathematics Studies~82.
Princeton University Press, 1975.

\smallskip
\item {[3]}
G.G. Gurzo.
{\it The group of smooth braids.}
16th All-Union Algebra Conference, Abstract II, 
39--40, Leningrad 1981.

\smallskip
\item {[4]}
Magnus, Karrass, and Solitar. 
``Combinatorial Group Theory.''
Dover Publications, New York, 1976.

\smallskip
\item {[5]}
D.L. Johnson.
{\it Towards a characterization of smooth braids.}
Mathematical Proceedings of the
Cambridge Philosophical Society 92 (1982) 425--427.

\smallskip
\item {[6]}
``The Kourovka notebook of unsolved problems in group
theory.''  7th edition, Novosibirsk, 1980.

\smallskip
\item {[7]}
H. W. Levinson.
{\it Decomposable braids and linkages}.
Transactions of the American Mathematical Society 178 
(1973), 111--126.

\smallskip
\item {[8]}
H. W. Levinson.
{\it Decomposable braids as subgroups of braid groups}.
Transactions of the American Mathematical Society 202 (1975), 51--55.

\smallskip
\item {[4]}
Magnus, Karrass, and Solitar.
``Combinatorial group theory. 
Presentations of groups in terms of generators and relations.''
Second revised edition.
Dover Publications, New York, 1976.

\smallskip
\item {[10]}
K. Y. Ng and T. Stanford.
{\it On Gusarov's groups of knots.}
Mathematical Proceedings of the 
Cambridge Philosophical Society 126 (1999), 63--76.

\smallskip
\item {[11]}
D. Rolfsen.
``Knots and Links''
Mathematics Lecture Series 7.
Publish or Perish, Wilmington, DE, 1976.

\end